\documentclass[11pt]{article}
\usepackage{latexsym}
\title{Spherical Designs and Generalized Sum-Free Sets in Abelian Groups}

\author{B\'{e}la Bajnok\\ Department of Mathematics\\ Gettysburg College\\ Gettysburg, PA 17325-1486 USA \\bbajnok@gettysburg.edu}

\newtheorem{thm}{Theorem}
\newtheorem{defin}[thm]{Definition}
\newtheorem{lem}[thm]{Lemma}
\newtheorem{cor}[thm]{Corollary}
\newtheorem{prop}[thm]{Proposition}

\newcommand{\bc}[2]{{{#1}\choose{#2}}}

\begin{document}

\renewcommand{\today}{September 15, 1999}

\maketitle

\begin{abstract}

We extend the concepts of sum-free sets and Sidon-sets of combinatorial number theory with the aim to provide explicit constructions for spherical designs. We call a subset $S$ of the (additive) abelian group $G$ {\it $t$-free}  if for all non-negative integers $k$ and $l$ with $k+l \leq t$, the sum of $k$ (not necessarily distinct) elements of $S$ does not equal the sum of $l$ (not necessarily distinct) elements of $S$ unless $k=l$ and the two sums contain the same terms. 

Here we shall give asymptotic bounds for the size of a largest $t$-free set in ${\bf Z}_n$, and for $t \leq 3$ discuss how $t$-free sets in ${\bf Z}_n$ can be used to construct spherical $t$-designs.

\end{abstract}

\section{Introduction}

In the attempt to provide explicit constructions for spherical designs, we introduce the concept of $t$-free sets (generalized sum-free sets) in abelian groups. $t$-free sets give an extention of the well studied concepts of sum-free sets, Sidon-sets, and $B_h$ sequences, and are the sources of some interesting number theory.

Section 2 of this paper gives a brief introduction to spherical designs, Section 3 describes their connection to $t$-free sets, and Section 4 gives some results on $t$-free sets. Readers only interested in $t$-free sets may proceed directly to the self-contained Section 4.

\section{Spherical Designs}

Spherical designs were introduced by Delsarte, Goethals, and Seidel in 1977 \cite{DelGoeSei:1977a}. 

\begin{defin} \label{designs}

A finite set $X$ of points on the $d$-sphere $S^d$ is a \emph{spherical $t$-design} or a \emph{spherical design of strength $t$}, if for every polynomial $f$ of total degree $t$ or less, the average value of $f$ over the whole sphere is equal to the arithmetic average of its values on $X$. If this only holds for 
homogeneous polynomials of degree $t$, then $X$ is called a \emph{spherical design of index $t$}.

\end{defin}

In other words, $X$ is of index $t$ if the Chebyshev-type quadrature formula
\begin{equation} \label{eq:quad}
\frac{1}{\sigma_d(S^d)} \int_{S^d}f({\bf x}) d \sigma_d({\bf x}) \approx \frac{1}{|X|} \sum_{{\bf x} \in X} f({\bf x})
\end{equation}
is exact for all homogeneous polynomials $f({\bf x})=f(x_0,x_1,\dots,x_d)$ of degree $t$ ($\sigma_d$ denotes the surface measure on $S^d$). $X$ is a $t$-design if it is of index $k$ for every $k \leq t$.  

The concept of $t$-designs on the sphere is analogous to $t-(v,k,\lambda)$ designs \cite{Sei:1990a}, and has been studied in various contexts, including representation theory, combinatorics, and approximation theory. For general references see \cite{Ban:1988a}, \cite{DelGoeSei:1977a}, \cite{God:1993a}, \cite{GoeSei:1979a}, \cite{GoeSei:1981a}, \cite{Hog:1996a}, \cite{Rez:1992a}, and \cite{Sei:1996a}. The existence of spherical designs for every $t$ and $d$ and large enough $n=|X|$ was first proved by Seymour and Zaslavsky in 1984 \cite{SeyZas:1984a}.

A central question in the field is to find all integer triples $(t,d,n)$ for which a spherical $t$-design on $S^d$ exists consisting of $n$ points, and to provide explicit constructions for these parameters. Delsarte, Goethals, and Seidel \cite{DelGoeSei:1977a} provide the tight lower bound
\begin{equation}
n \geq \bc{d+\lfloor t/2 \rfloor}{\lfloor t/2 \rfloor}+\bc{d+\lfloor (t-1)/2 \rfloor}{\lfloor (t-1)/2 \rfloor}=O(t^d). \label{eq:lower}
\end{equation}
We shall refer to the bound \ref{eq:lower} as the DGS bound. Spherical designs of this minimum size are called \emph{tight}. Bannai and Damerell \cite{BanDam:1979a}, \cite{BanDam:1980a} proved that tight spherical designs for $d \geq 2$ exist only for $t=1,2,3,4,5,7$, or $11$. All tight $t$-designs are known, except possibly for $t=4,5,$ or 7. In particular, there is a unique $11-$design ($d$=23 and $n=196,560$).

The first general construction of spherical designs for arbitrary $t$, $d$, and large enough $n$ were given independently by Wagner \cite{Wag:1991a} and the author \cite{Baj:1991d}, \cite{Baj:1992a}, who used $n \geq C(d)t^{O(d^4)}$ and $n \geq C(d)t^{O(d^3)}$ points, respectively. This bound was later reduced to $C(d)t^{d^2/2+d/2}$ by Korevaar and Meyers \cite{KorMey:1994a}. They believe that the minimum size of a $t$-design on $S^d$ is $C(d)t^d$.

It is easy to see that a 1-design of size $n$ exists on $S^d$ for every $d$ and $n \geq 2$ (take any point set whose centroid is the origin). The case $t=2$ was settled by Mimura \cite{Mim:1990a} who proved the following.

\begin{thm} \label{Mimura}

A 2-design of size $n$ exists on $S^d$ if and only if $n \geq d+2$ (the DGS bound (\ref{eq:lower})), unless $n$ is odd and $n=d+3$. 

\end{thm}

For $t=3$ the author \cite{Baj:1998a} provided constructions for all $d$ and $n$ for which 3-designs are believed to exist. Namely, we have the following.

\begin{thm} \label{t=3} 

A 3-design of size $n$ on $S^d$ exists for every $n \geq 2(d+1)$ (the DGS bound (\ref{eq:lower})), unless $n$ is odd and $n<5(d+1)/2$ or $(d,n)\in \{(2,9),(4,13)\}$. 

\end{thm}

We conjectured in \cite{Baj:1998a} that 3-designs do not exist for other parameters. This conjecture is supported by the numerical evidence of Hardin and Sloane \cite{HarSlo:1996a} and by a result of Boyvalenkov, Danyo, and Nikova \cite{BoyDanNik:1999a} that no 3-design exists of size $n$ on $S^d$ if $n$ is odd and $n<(\sqrt[3]{2}+1)(d+1)+.3 \approx 2.26(d+1)+.3$. In particular, there is no 3-design on 7 points on $S^2$, leaving only the case $n=9$ open on the 2-dimensional sphere. 

The proof of Theorem \ref{t=3} is based on a number theoretic idea. Below we shall describe this method in a general setting that might be of independent interest.

\section{From Spherical Designs to Additive Number Theory} 

For explicit constructions of spherical designs it is convenient to use the following equivalent definition, see \cite{DelGoeSei:1977a} or \cite{Ban:1988a}.

\begin{lem} \label{equiv}

A finite subset $X$ of $S^d$ is a spherical $t$-design if and only if for every homogeneous harmonic polynomial $f$ of total degree $t$ or less $$\sum_{{\bf x} \in X} f({\bf x})=0.$$

\end{lem}

A polynomial $f(x_0,x_1,\dots,x_d)$ is \emph{harmonic} if it satisfies Laplace's equation $\Delta f=0$. The set of homogeneous harmonic polynomials of degree $k$ over $S^d$ forms the vector space $Harm_{k}(S^d)$ with $$\dim Harm_{k}(S^d)=\bc{d+k}{k}-\bc{d+k-2}{k-2}.$$ In particular, for $k=1,2$, and $3$ we find that $\Phi_k(S^d)$ forms a basis for $Harm_k(S^d)$ where

$$\Phi_1(S^d)=\{x_i | 0 \leq i \leq d\},$$
$$\Phi_2(S^d)=\{x_ix_j| 0 \leq i<j \leq d\} \cup \{x_i^2-x_{i+1}^2 | 0 \leq i \leq d-1\},$$ and
$$\Phi_3(S^d)=\{x_ix_jx_k | 0 \leq i<j<k \leq d\} \cup \{x_i^3-3x_ix_j^2 | 0 \leq i \not =j\leq d\}.$$

We now attempt to find a set of $n$ points on $S^d$ which forms a $t$-design. Before we proceed, we state the following lemma. 

\begin{lem} \label{lem}

For all positive integers $a$ and $n$ we have $\sum_{i=1}^n \sin(\frac{2 \pi i}{n}a)=0$. Furthermore, if $a$ is not a multiple of $n$, then $\sum_{i=1}^n \cos(\frac{2 \pi i}{n}a)=0$ as well. 

\end{lem}

{\it Proof.} Let $z=\cos(\frac{2 \pi}{n}a)+\sqrt{-1} \sin (\frac{2 \pi}{n}a)$, the $a$-th complex value of $\sqrt[n]{1}$. If $a$ is not a multiple of $n$ then $z \not =1$, and we have $\sum_{i=1}^n z^i=0$. $ \quad \Box$

For $t=1$ the lower bound (\ref{eq:lower}) yields $n \geq 2$. By Lemma \ref{lem} we see that the vertices $${\bf u_i}=\left(\cos(\frac{2 \pi i}{n}), \sin (\frac{2 \pi i}{n}),0,\dots,0\right)$$ ($i=1,2,\dots,n$) of a regular $n$-gon on the equator of $S^d$ form a $1$-design if $n\geq2$.

Below we shall try to generalize this simple construction to the case of $t \geq 2$. We follow methods similar to those used by Mimura \cite{Mim:1990a} and the author \cite{Baj:1998a}. For simplicity we assume in what follows that $d$ is odd and let $d=2m-1$. The case when $d$ is even can be reduced to this case by a simple technique, see \cite{Baj:1998a} or \cite{Mim:1990a}. 

Suppose that $a_1,a_2,\dots,a_m$ are positive integers, and consider the $n$ points $${\bf u_i}=\frac{1}{\sqrt{m}}\left(\cos(\frac{2 \pi ia_1}{n}), \sin (\frac{2 \pi ia_1}{n}),\dots,\cos(\frac{2 \pi ia_m}{n}), \sin (\frac{2 \pi ia_m}{n})\right)$$
where $i=1,2,\dots,n$. The scalar in front is chosen so that each point is on $S^d$.

We now examine $\sum_{i=1}^n f({\bf u_i})$ for monomials $f$ of $d+1=2m$ variables. Suppose that $f=\prod_{j=0}^dx_j^{k_j}$ with $1 \leq \sum_{j=0}^d k_j=k$. Then 
$$f({\bf u_i})=\frac{1}{m^{k/2}}\cos^{k_0}(\frac{2 \pi ia_1}{n})\sin^{k_1}(\frac{2 \pi ia_1}{n})\cdots\cos^{k_{d-1}}(\frac{2 \pi ia_m}{n})\sin^{k_d}(\frac{2 \pi ia_m}{n}).$$

Using the trigonometric identities
$$\sin \alpha \sin \beta = \frac{1}{2}[\cos(\alpha - \beta) - \cos(\alpha + \beta)],$$
$$\sin \alpha \cos \beta = \frac{1}{2}[\sin(\alpha - \beta) + \sin(\alpha + \beta)],$$ and
$$\cos \alpha \cos \beta = \frac{1}{2}[\cos(\alpha - \beta) + \cos(\alpha + \beta)]$$
 repeatedly, we can write $f({\bf u_i})$ as the sum of $2^{k-1}$ terms, each of the form $$\pm \frac{1}{2^{k-1}m^{k/2}} \sin \left( \frac{2 \pi i}{n}(l_1a_1 +\cdots + l_ma_m) \right)$$ or $$\pm \frac{1}{2^{k-1}m^{k/2}} \cos \left( \frac{2 \pi i}{n}(l_1a_1 + \cdots + l_ma_m) \right),$$ where $l_1,l_2,\dots,l_m$ are integers with $|l_{\nu}| \leq k_{2{\nu}-2}+k_{2{\nu}-1}$ for $\nu=1,2,\dots,m$; in particular, $|l_1|+ \cdots + |l_m| \leq k$. In fact, a closer look reveals that if either $k_{2{\nu}-2}$ or $k_{2{\nu}-1}$ is odd, then it is possible to do this so that a cosine term with $l_{\nu}=0$ does not appear; in particular, a cosine term with $l_1=l_2= \cdots = l_m=0$ will not appear if at least one exponent $k_j$ is odd ($j=0,1, \dots,d$).

Therefore, by Lemma \ref{lem}, we have the following theorem.

\begin{thm} \label{main}

Let $f$ be a monomial of $2m$ variables, total degree $k$, and suppose that at least one variable appears with an odd exponent in $f$. Suppose further that the integers $a_1,a_2,\dots,a_m$ are chosen so that $l_1a_1 + \cdots + l_ma_m$ is not divisible by $n$ whenever the integers $l_1,l_2,\dots,l_m$ satisfy $1 \leq |l_1|+ \cdots + |l_m| \leq k$. For $i=1,2,\dots,n$ define $${\bf u_i}=\frac{1}{\sqrt{m}}\left(\cos(\frac{2 \pi ia_1}{n}), \sin (\frac{2 \pi ia_1}{n}),\dots,\cos(\frac{2 \pi ia_m}{n}), \sin (\frac{2 \pi ia_m}{n})\right).$$ Then we have $\sum_{i=1}^n f({\bf u_i})=0$. $ \quad \Box$

\end{thm}

A set $\{a_1,a_2,\dots,a_m\}$ of integers satisfying the condition in Theorem \ref{main} will be called \emph{$k$-free}. This concept leads us to the beautiful area of additive number theory, and shall be discussed in the next Section.

\begin{cor} \label{cor}

If $k$ is an odd positive integer and the set $\{a_1,a_2,\dots,a_m\}$ is $k$-free, then $X=\{{\bf u_i}|i=1,2,\dots,n\}$ (as defined in Theorem \ref{main} above) is a spherical design on $S^d$ of index $k$. $ \quad \Box$

\end{cor}

Finally, as in \cite{Baj:1998a}, we find that for $f=x_j^2$ ($j=0,1, \dots,d$) and a 2-free set $\{a_1,a_2,\dots,a_m\}$, $\sum_{i=1}^n f({\bf u_i})=n/2$. Since this value is independent of $j$, we see that $\sum_{i=1}^n f({\bf u_i})=0$ for every polynomial in $\Phi_1(S^d) \cup \Phi_2(S^d) \cup \Phi_3(S^d)$, thus we get the following. 

\begin{cor} \label{this}

Let $t=1,2$, or 3. Suppose that the set $\{a_1,a_2,\dots,a_m\}$ is $t$-free. Then $X=\{{\bf u_i}|i=1,2,\dots,n\}$ (as defined in Theorem \ref{main} above) is a spherical $t$-design on $S^d$. $ \quad \Box$

\end{cor}

The earlier stated Theorem \ref{t=3} is based on Corollary \ref{this}. The application of our methods to $t \geq 4$ will be the subject of further study.

\section{Generalized Sum-Free Sets in Abelian Groups}

In this section $t$ is a positive integer, and $G$ is an abelian group written in additive notation. In view of Theorem \ref{main}, we make the following definition. 

\begin{defin} \label{t-free}

We say that $S \subset G$ is a \emph{$t$-free set} in $G$ if for all non-negative integers $k$ and $l$ with $k+l \leq t$, the sum of $k$ (not necessarily distinct) elements of $S$ can only equal the sum of $l$ (not necessarily distinct) elements of $S$ if $k=l$ and the two sums contain the same terms.

\end{defin}

Equivalently, we say that $S$ is $t$-free in $G$ if every equation of the form $\epsilon_1 x_1+\epsilon_2 x_2+\cdots+\epsilon_t x_t=0$, where $\epsilon_i=0,\pm1$ for $i=1,2,\dots,t$, has only trivial solutions in $S$: each $\epsilon_i=0$, or the same $x_i$ appears both with a coefficient of $1$ and of $-1$. The cardinality of a largest $t$-free set in $G$ will be denoted by $s(G,t)$.

Our $t$-free sets are extensions of the extensively studied concepts of sum-free sets and Sidon sets in abelian groups. A {\it sum-free set in $G$} is a subset $S$ of $G$ for which $(S +S) \cap S= \emptyset$, i.e. there are no (not necessarily distinct) elements $a$, $b$, and $c$ in $S$ for which $a+b=c$. A {\it Sidon set in $G$} is a subset $S$ of $G$ for which the only way to have (not necessarily distinct) $a,b,c,d \in S$ with $a+b=c+d$ is the trivial $\{a,b\}=\{c,d\}$. Sum-free sets, Sidon sets, and their generalizations such as $B_h$ sequences have a long history and have been investigated extensively, most notably by Erd\H{o}s. For more information see, for example, \cite{AloKle:1990a}, \cite{ErdFre:1991a}, \cite{Guy:1994a}, \cite{HalRot:1983a}, \cite{WalStrWal:1972a}, and their references.  

Here we are interested in $t$-free sets in the group ${\bf Z}_n$. For an explicit construction and to have a lower bound on $s({\bf Z}_n,t)$, we have the following.

\begin{prop} \label{log} If $n>t^m$, then the set $\{1,t,t^2,\dots,t^{m-1}\}$ is a $t$-free set of size $m$ in ${\bf Z}_n$. This gives $s({\bf Z}_n,t) \geq \lfloor \log_t (n-1) \rfloor$. $ \quad \Box$

\end{prop}

We can certainly find better approximations for $s({\bf Z}_n,t)$. For $t=1$ we can obviously take the set $\{1,2,\dots,n-1\}$, hence $s({\bf Z}_n,1)=n-1$. It is also clear that $\{1,2,\dots, \lfloor (n-1)/2 \rfloor \}$ is a 2-free set and, since we can never have both $a$ and $n-a$ in a 2-free set, we conclude that $s({\bf Z}_n,2)=\lfloor (n-1)/2 \rfloor$. 

For $t=3$ the situation becomes interesting. First we prove the following.

\begin{prop} \label{n/4}

For every $n$ we have $s({\bf Z}_n,3) \leq \lfloor \frac{n}{4} \rfloor$.

\end{prop}

{\it Proof.} [Based on \cite{Boy}.] As above, we note that we can assume without loss of generality that our 3-free set $S=\{a_1,a_2,\dots,a_m\}$ is such that $1 \leq a_1 < a_2 < \cdots < a_m \leq \lfloor \frac{n-1}{2} \rfloor$. Consider the set $S^*=\{a_1,a_2,\dots,a_m,a_m-a_1,a_m-a_2,\dots,a_m-a_{m-1}\} \subseteq \{1,2,\dots,\lfloor \frac{n-1}{2} \rfloor\}$. Since $S$ is 3-free, the $2m-1$ elements in $S^*$ are all distinct in ${\bf Z}_n$. Our claim then holds if $|S^*| \leq \lfloor \frac{n}{2} \rfloor -1$. This is clearly the case if $n$ is even.

When $n$ is odd, we argue as follows. If $S^*$ were equal to $\{1,2,\dots,\lfloor \frac{n-1}{2} \rfloor\}$, then we would have to have $a_m=\frac{n-1}{2}$. If $a_1=1$, then $a_1+a_m+a_m=n$, a contradiction, so assume that $a_1 > 1$. Now $1,2,\dots, a_1-1 \in S^*$ can only happen if $a_1 \leq m$ and $1=a_m-a_{m-1}, 2=a_m-a_{m-2}, \dots, a_1-1=a_m-a_{m-(a_1-1)}$. But then the three elements $a_m, a_{m-(a_1-1)}$, and $a_1$ add to $2a_m+1=n$, a contradiction again. $ \quad \Box$

We indeed have a 3-free set in ${\bf Z}_n$ of size $\lfloor \frac{n}{4} \rfloor$ when $n$ is even: take the odd integers up to (but not including) $n/2$ (note that three odd numbers will not add to $n$ if $n$ is even). When $n$ is odd, we can still take the set of odd integers up to (but not including) $n/3$. Surprisingly, we can do better in one case, namely if $n$ has a divisor $p$ of the form $p=6q+5$. In this case, the set $$ \{ ip+2j+1 | i=0,1, \dots, \frac{n}{p}-1, j=0,1, \dots,q \}$$ is 3-free (see \cite{Baj:1998a}). Thus we have the following proposition.

\begin{prop} \label{3-free}

If $n$ is even, then $s({\bf Z}_n,3) = \lfloor \frac{n}{4} \rfloor$.
If $n$ is odd and has no prime divisors congruent to 5 (mod 6), then $s({\bf Z}_n,3) \geq \lfloor \frac{n}{6} \rfloor$.
If $n$ is odd and $p$ is its smallest prime divisor congruent to 5 (mod 6), then $s({\bf Z}_n,3) \geq \frac{p+1}{6p}n$. $ \quad \Box$

\end{prop}

We mention here that by a recent result of Ruzsa \cite{Ruz:1999a} the lower bounds for $s({\bf Z}_n,3)$ given in Theorem \ref{3-free} are exact for every $n$.

Note that by Propositions \ref{n/4} and \ref{3-free}, $\lfloor \frac{n}{6} \rfloor \leq s({\bf Z}_n,3) \leq \lfloor \frac{n}{4} \rfloor$. For $t \geq 4$, however, we have $s({\bf Z}_n,t)=o(n)$. Below we establish the following asymptotic results.

\begin{thm} \label{asymp}

For a given positive integer $t$, there are constants $c_1(t)$ and $c_2(t)$ for which $c_1(t) n^{1/t} \leq s({\bf Z}_n,t) \leq c_2(t) n^{1/\lfloor t/2 \rfloor}$ for every positive integer $n$.

\end{thm}

The lower and upper bounds in Theorem \ref{asymp} will follow from Propositions \ref{lower} and \ref{upper}, respectively. 

\begin{prop} \label{lower}

Let $t$, $m$, and $n$ be positive integers for which $n \geq t3^tm^t$. Then ${\bf Z}_n$ has a $t$-free set of size $m$.

\end{prop}

\begin{prop} \label{upper}

Let $t$, $m$, and $n$ be positive integers for which ${\bf Z}_n$ has a $t$-free set of size $m$. Then $n \geq \bc{m+\lfloor t/2 \rfloor}{\lfloor t/2 \rfloor}$.

\end{prop}

{\it Proof of Proposition \ref{lower}.} We use induction on $m$. For $m=1$ we see that $\{1\}$ is a $t$-free set in ${\bf Z}_n$ whenever $n>t$, and this indeed holds by our assumption.

Assume now that our proposition holds for a positive integer $m$ and suppose that $n \geq t3^t(m+1)^t$. Since this value is greater than $t3^tm^t$, our inductive hypothesis implies that ${\bf Z}_n$ has a $t$-free set $S$ of size $m$.

Let us define $\Gamma^t S:=\{\epsilon_1 s_1+\epsilon_2 s_2+\cdots+\epsilon_t s_t=0 | \epsilon_i \in \{0,1,-1\}, s_i \in S, i=1,2,\dots,t \}$. We have $|\Gamma^t S| \leq 3^tm^t$.

Now look at $A_j:=\{j,2j,\dots,tj\} \subset {\bf Z}_n$ for $j=1,2,\dots,3^tm^t+1$. Since $|\Gamma^t S| \leq 3^tm^t$, we must have a $j_0 \in \{1,2,\dots,3^tm^t+1\}$ for which $\Gamma^t S \cap A_{j_0} = \emptyset$. We claim that $S \cup \{j_0\}$ is a $t$-free set in ${\bf Z}_n$ of size $m+1$.

First, $|S \cup \{j_0\}|=m+1$, since $j_0 \in S$ implies $j_0 \in \Gamma^t S \cap A_{j_0}$, a contradiction. To show that $S \cup \{j_0\}$ is $t$-free, assume that $x+kj_0=0$ for some $x \in \Gamma^t S$ and $0 \leq k \leq t$. If $x \not =0$, then $k \not =0$ and $-x=kj_0 \in A_{j_0}$. But $-x \in \Gamma^t S$ (as $\Gamma^tS$ is closed under taking negatives), hence $-x \in \Gamma^t S \cap A_{j_0} = \emptyset$, a contradiction. On the other hand, if $x=0$, then $kj_0=0$. But $x=0$ implies that $x$ is trivial (because $S$ is $t$-free), so we just need to prove that $k=0$ as well. This indeed holds, as $1 \leq k \leq t$ would imply that the positive integer $kj_0$ is at most $t(3^tm^t+1)$, a number less than $t3^t(m+1)^t$, so $kj_0 \not =0$ in ${\bf Z}_n$, a contradiction. $ \quad \Box$

{\it Proof of Proposition \ref{upper}.} Let $S$ be a $t$-free set in ${\bf Z}_n$. Consider the set $\Sigma^{\lfloor t/2 \rfloor} S=\cup_{k=1}^{\lfloor t/2 \rfloor}  \underbrace{S+S+\cdots+S}_k$. 

Since $S$ is $t$-free, $a_1+a_2+\cdots +a_i \in \Sigma^{\lfloor t/2 \rfloor} S$ and $b_1+b_2+\cdots +b_j \in \Sigma^{\lfloor t/2 \rfloor} S$ are different and non-zero in ${\bf Z}_n$, unless $i=j$ and $\{a_1,a_2,\dots,a_i\}=\{b_1,b_2,\dots,b_j\}$ as multisets. Therefore, $n-1 \geq |\Sigma^{\lfloor t/2 \rfloor} S|=\sum_{k=1}^{\lfloor t/2 \rfloor} \bc{m+k-1}{k} = \bc{m+\lfloor t/2 \rfloor}{\lfloor t/2 \rfloor}-1$, from which our claim follows. $ \quad \Box$

We close our paper with some open problems on $t$-free sets.

{\bf Problem 1.} Find the correct asymptotic value of $s({\bf Z}_n,t)$ (see Theorem \ref{asymp}).

{\bf Problem 2.} Improve Proposition \ref{log} by finding an explicit construction for a $t$-free set of size $m$ in ${\bf Z_n}$ for $m$ at least as in Proposition \ref{lower}.

{\bf Problem 3.} Investigate $t$-free sets in other abelian groups and (after a modified definition) in non-abelian groups.

\end{document}